\documentclass[12pt,a4paper]{article}
\usepackage[T2A]{fontenc}
\usepackage[cp1251]{inputenc}
\usepackage[russian, english]{babel}
\usepackage{indentfirst}
\usepackage{amsmath,amssymb,amsthm,amsfonts,amscd}
\usepackage{color}

\newtheorem{teo}{Theorem}

\newtheorem{lem}{Lemma}

\newtheorem{prp}{Proposition}

\newtheorem{dfn}{Definition}

\newtheorem{rem}{Remark}

\title{On hyperbolic attractors and repellers of endomorphisms}
\author{V.~Z.~Grines\footnote{
Viacheslav\,Z.~Grines;
E-mail: vgrines@hse.ru;
National Research University Higher School of Economics,
Nizhniy Novgorod, Russia.},
E.~D.~Kurenkov\footnote{
Evgeniy\,D.~Kurenkov
E-mail: ekurenkov@hse.ru;
National Research University Higher School of Economics,
Nizhniy Novgorod, Russia.}}
\date{}
\begin{document}

\maketitle

\vskip 0.25 true cm

\noindent \textbf{Abstract.} It is well known that topological classification of dynamical systems with hyperbolic dynamics is significantly defined by dynamics on nonwandering set.  F.~Przytycki generalized axiom $A$ for smooth endomorphisms that was previously introduced by S.~Smale for diffeomorphisms and proved spectral decomposition theorem which claims that nonwandering set of an $A$-endomorphism is a union of a finite number basic sets.
In present paper the criterion for a basic sets of an $A$-endomorphism to be an attractor is given. 
Moreover, dynamics on basic sets of codimension one is studied. It is shown, that if an attractor is a topological submanifold of codimension one of type $(n-1, 1)$, then it is smoothly embedded in ambient manifold and restriction of the endomorphism to this basic set is an expanding endomorphism. If a basic set of type $(n, 0)$ is a topological submanifold of codimension one, then it is a repeller and restriction of the endomorphism to this basic set is also an expanding endomorphism.

\vskip 0.5 true cm

\noindent \textbf{Keywords}: endomorphism, axiom $A$, basic set, attractor, repeller.

\section{Introduction}

It is well known that axiom $A$ introduced by S.~Smale along with  strong transversality condition are necessary and sufficient for structural stability of a dynamical system (either smooth flow or diffeomorphism) given on a smooth manifold. There are several classification results for such systems based on Smale's theorem on  spectral decomposition which states that the nonwandering set of a structurally stable system can be uniquely decomposed into finitely many closed invariant basic sets each of which contains  a transitive trajectory (see~\cite{An:Grub,GrZh,GrPo,GrMePoLe}  for information and references).

As for noninvertible discrete dynamical systems (endomorphisms), there are only few classes of systems satisfying axiom $A$ with a well studied structure of basic sets. Among them there are endomorphisms of the circle and the interval~\cite{Van,Jak}, endomorphisms of the Riemann sphere that occur in holomorphic dynamics~\cite{Mi,Lu}, expanding endomorphisms of closed manifolds~\cite{Shu}. In the present paper we study dynamics of endomorphisms on basic sets of codimension one which are topological submanifolds of ambient manifold.

Let $M^n$ be a smooth closed manifold. By $C^r$-endomorphism we mean a $C^r$-smooth, $r \geqslant 1$, surjective map   $f\colon M^n \to M^n$. If an endomorphism $f$ has a  $C^r$-smooth inverse, then $f$ is called a $C^r$-diffeomorphism.

Let $f \colon M^n \to M^n$ be an $C^r$-endomorphism. Let us define $\hat M$ the subset of Tikhonov product $ \widetilde M=\prod \limits_ {i = -\infty}^{+ \infty} M^n $ as $\hat M = \left\{ \{x_i\}_{i \in \mathbb Z} \in \widetilde M \mid f(x_i) = x_{i+1} \right\}$. For any $x\in M^n$ set   $\hat x = \left\{ \{x_i\}_{i \in \mathbb Z} \in \hat M \mid x_0 = x \right\}$.  

If $f$ is a diffeomorphism, then for any point $x \in M^n$ the set  $\hat x$ consists of exactly one point. If $f$ is not one-to-one, then it is not true in general.
Let  use symbol $\bar x$ for  a particular element of $\hat x$.
For an $f$-invariant set $\Lambda$ (i. e. $f(\Lambda) = \Lambda$) let us introduce $\hat \Lambda \subset \hat M$ as $\hat \Lambda = \left\{ \{x_i\}_{i \in \mathbb Z} \in \hat M \mid x_i \in \Lambda, \, \forall i \in \mathbb Z \right\}$. For $\bar x \in \hat M$ let $\hat f(\bar x)$ be a shift of $\bar x$ (i.e. $\hat f(\{x_i\}_{i \in \mathbb Z}) = \{x_{i+1}\}_{i \in \mathbb Z}$; $\hat f^{-1}(\{x_i\}_{i \in \mathbb Z}) = \{x_{i-1}\}_{i \in \mathbb Z}$).
Let $\Lambda\subset M^n$ be closed $f$-invariant set. The following definition of hyperbolic set given by F.~Przytycki~\cite{Prz1} generalizes Smale's definition for diffeomorphisms~\cite{Sm}.

\begin{dfn} Let $f$ be an endomorphism of a manifold $M^n$. An invariant set $\Lambda$ is called hyperbolic if there exist constants $C > 0$, $0<\lambda<1$ such that for every $x \in \Lambda$ and  every $\bar x \in \hat x \cap \hat \Lambda$ ($\bar x  = \{x_i\}_{i \in Z}$)  there exists a continuous splitting of the tangent subbundle $\bigcup\limits_{i \in \mathbb Z}T_{x_i}M^n$ into the direct sum $\bigcup\limits_{i \in \mathbb Z}T_{x_i}M^n = \bigcup\limits_{i \in \mathbb Z} E^s_{x_i, \hat f^i(\bar x)} \oplus E^u_{x_i, \hat f^i(\bar x)}$ such that:
\begin{enumerate}
\item $Df\left(E^s_{x_i, \hat f^i(\bar x)}\right) = E^s_{x_{i+1}, \hat f^{i+1}(\bar x)}$, $Df\left(E^u_{x_i, \hat f^i(\bar x)}\right) = E^u_{x_{i+1}, \hat f^{i+1}(\bar x)}$, where $E^s_{x_i, \hat f^i(\bar x)}, E^u_{x_i, \hat f^i(\bar x)} \subset T_{x_i}M^n$;
\item $\|Df^k(v)\| \leqslant C\lambda^k\|v\|$, for all $k \geqslant 0, i \in \mathbb Z$, $v \in E^s_{x_i, \hat f^i(\bar x)}$;
\item$\|Df^k(v)\| \geqslant(1/C)\lambda^{-k}\|v\|$, for all $k \geqslant 0, i \in \mathbb Z$, $v \in E^u_{x_i, \hat f^i(\bar x)}$.
\end{enumerate}
\end{dfn}

\begin{rem}
One can show that the  stable subspace $E^s_{x_i, \hat f^i(\bar x)}$ in the tangent space $T_{x_i}M^n$ at point $x_i$ is independent of  the choice of   $\bar x \in \hat x$ (see, for example,~\cite{Prz1}). Also note, that $E^s_{x_0, \bar x} = E^s_{x_0, \hat f^{0}(\bar x)}$ and $E^u_{x_0, \bar x} = E^u_{x_0, \hat f^{0}(\bar x)}$.
\end{rem}

For  a smooth map $f\colon M^n\to M^n$,     a point  $x\in M^n$ is called {\it regular}   if the rank of the map $Df(x) : T_x M^n \to T_{f(x)}M^n$ is exactly $n$. Otherwise  the point $x$ is called {\it 
singular}.

Next defintion is a generalization of Smale's axiom $A$ that was given in the paper \cite{Prz1}.

\begin{dfn}\label{axiom_a}
An endomorphism $f \colon M^n \to M^n$  satisfies  axiom $A$ if the following conditions hold:
\begin{enumerate}
\item the nonwandering set $\Omega_f$ is hyperbolic and does not contain singular points of $f$;
\item the set of periodic points $Per_f$ of $f$ is dense in nonwandering set $\Omega_f$.
\end{enumerate}
\end{dfn}

For $A$-endomorphisms there exists a spectral decomposition theorem proved in~\cite{Prz1} that generalizes Smale's result for diffeomorphisms~\cite{Sm}.

\begin{prp}\label{spectral}
Let $f$ be an $A$-endomorphism. Then its nonwandering set 
$\Omega_f$ can be uniquely decomposed into a finite union of closed $f$-invariant sets (called basic sets) $\Omega_f = \bigcup\limits_{i = 1}^l \Lambda_i$ such that the restriction of $f$ to every basic set $\Lambda_i$ is topologically transitive. 
\end{prp}

\begin{dfn}
A basic set $\Lambda$ of an  $A$-endomorphism $f$ is called  attractor if it has a closed neighborhood $U \supset \Lambda$ such that $f(U)\subset \mathop{\rm{Int}} U$ and $\bigcap\limits_{n = 0}^\infty f^n(U) = \Lambda$.
\end{dfn}

\begin{dfn}\label{def_repeller}
A basic set $\Lambda$ of an $A$-endomorphism $f$ is called  repeller if  it has an open neighborhood $U$ such that $\mathop{\text{cl}\,(U)} \subset f(U)$ and  $\bigcap\limits_{n = 0}^\infty f^{-n}(U) = \Lambda$.
\end{dfn}

The following definition belongs to M.~Shub~\cite{Shu}.

\begin{dfn}\label{Shu_expand}
A $C^r$-endomorphisms $f \colon M^n \to M^n$ is called expanding if there exist constants $C > 0$ and $\lambda >1$ such that $\|Df^n(v)\| \geqslant C\lambda^n \|v\|$ for all $v \in TM^n$, $n = 0, 1, 2, \ldots$.
\end{dfn}
   
Nevertheless, is is possible to define an expanding endomorphism not only for smooth manifolds but for arbitrary metric space as well. The following definition was given in~\cite{Kat}.

\begin{dfn}\label{Kat_expand}
A continuous map $f \colon X \to X$ of a metric space $X$  is called expanding if there exist constants $\varepsilon > 0$ and $\mu > 1$ such that for all $x, y \in X$, $x \neq y$, $\rho(x, y) < \varepsilon$ the following inequality holds $\rho(f(x), f(y)) > \mu \rho(x, y)$.
\end{dfn}
Note that in the case when $X$ is a $C^1$-smooth,   compact manifold and   $f$ is a $C^1$-smooth map it follows from \cite{Kat} that conditions  of definition \ref{Shu_expand} imply  conditions of definition  \ref{Kat_expand}.

It follows from definition \ref{Shu_expand} that an ambient manifold $M^n$ of an expanding endomorphism $f$ is hyperbolic set. Moreover, it was shown in \cite{Shu} that if is $M^n$ is compact, then periodic points of expanding endomorphism $f$ are dense in $M^n$. Thus, any expanding endomorphism of a compact manifold is  $A$-endomorphism and its nonwandering set coincides with the  ambient manifold.  It was also shown in this  paper  that an ambient manifold of an expanding endomorphism has an Euler characteristics equal to zero and that its universal covering is diffeomorphic to $\mathbb R^n$. Moreover if a compact manifold $M^n$ is locally flat, then it admits  an expanding endomorphism \cite{EpShu},.

It was shown in \cite{Shu} that if ambient manifold $M^n$ is diffeomorphic to the   $n$-torus $\mathbb{T}^n$, then expanding endomorphism $f$ is topologically conjugated to  the  algebraic expanding endomorphism.

\begin{dfn}\label{anosendom}
An endomorphism $f \colon M^n \to M^n$ is called Anosov endomorphism if the  ambient manifold $M^n$ is hyperbolic set.
\end{dfn}

It follows from definition \ref{anosendom} that expanding endomorphism is an Anosov endomorphism.

Other examples of endomorphisms such that their unique basic set coincides with the  ambient manifold are provided by Anosov algebraic endomorphisms of $n$-torus induced by matrix $A_{n \times n}$ with eigenvalues which are   inside and outside of the unit circle and with no eigenvalues on the unit circle.

It is well known that an arbitrary Anosov diffeomorphism of $n$-torus is conjugate with the algebraic hyperbolic automorphism \cite{Fr,New,Mann}. However there is no such result for an Anosov endomorphism that are not expanding or diffemorphisms \cite{Prz1,MaPu}. Moreover it was shown in paper~\cite{Zhang}  that the set of endomorphisms  of the  $n$-torus that are not conjugated  to any algebraic endomorphism is the  residual\footnote{The set $X$ is called residual, if it is an intersection of countably many sets with dense interiors.} subset in the set of all  Anosov endomorphisms on the torus. Thus, the question whether an ambient manifold of an arbitrary Anosov endomorphism is nonwandering is still open.

As it follows from \cite{An,Shu} Anosov diffeomorphisms and expanding endomorphism are structurally stable. However, F.~Przytycki~\cite{Prz1} and R.~Mane with Ch.~Pugh \cite{MaPu}  independently proved that Anosov endomorphisms that are not expanding and diffeomorphisms are not structurally stable.

For a basic set $\Lambda$ of an $A$-endomorphism $f\colon M^n\to M^n$ the  pair of integers $(\dim E^u_{x_0, \bar x}, \dim E^s_{x_0, \bar x})$ is called {\it the  type of basic set}.  It follows from \cite{Prz1} that this definition is correct since  $\dim E^u_{x_0, \bar x}, \dim E^s_{x_0, \bar x}$ do not depend on the point $\bar x \in \hat \Lambda$.

In the case when $f$ is $A$-diffeomorphism the topological structure of basic sets of codimension one is well studied. It follows from R.~V.~Plykin that  any basic set  of codimension one is necessarily  either attractor or repeller. In the case when $n = 2$ it is locally homeomophic to the product of the Cantor set and the  interval. If $f \colon M^3 \to M^3$ is an $A$-diffeomorphism of closed $3$-manifold and $\Lambda$ is a 2-dimensional basic set of the type $(2, 1)$ ($(1, 2)$) that coincides with the union of unstable (stable) manifolds of its points, then it is called {\it expanding attractor (contracting repeller)}, and it follows from~\cite{Pl,Br}, that it is locally homeomorphic to the product of $\mathbb R^{2}$ and the Cantor set. If $\Lambda$ doesn't coincides with the union of unstable (stable) manifolds of its points, then $\Lambda$ is homeomorphic to the $2$-torus $\mathbb T^{2}$ and restriction of $f_\Lambda$ is topologically conjugated  to the Anosov algebraic diffeomorphism.

In the present paper we prove a  criterion for a basic set $\Lambda$  to be an attractor of an endomorphism $f \colon M^n \to M^n$. Moreover, we study the dynamics of restriction of an $A$-endomorphisms $f \colon M^n \to M^n$ to the basic set $\Lambda$ in the case when it is submanifold of codimension one of the manifold $M^n$. The existence of a smooth structure on such basic set is also considered.

\begin{teo}\label{criterium}
A basic set $\Lambda$  of an $A$-endomorphism $f \colon M^n\to M^n$ is an attractor if and only if there exists $\varepsilon > 0$ such  that for every $x \in \Lambda$ and  $\bar x \in \hat x \cap \hat \Lambda$   one has  $W^u_{x, \bar x, \varepsilon} \subset \Lambda$.\footnote{$W^u_{x, \bar x, \varepsilon}$ is introduced in section~\ref{aux} (definition~\ref{def:local}).}
\end{teo}

\begin{teo}\label{repeller}
Let basic set $\Lambda$ be a codimension one topological submanifold of $M^n$. If $\Lambda$ is of type $(n, 0)$, then:

1) $\Lambda$ is a repeller; 

2) the restriction of $f$ to $\Lambda$ is an expanding endomorphism.\footnote{In sense of definition~\ref{Kat_expand} in metric on $M^n$ induced by Riemannian metric $\langle \cdot, \cdot \rangle_{\Lambda}$ defined in proposition~\ref{metric}.} 
\end{teo} 

\begin{rem}
If in theorem \ref{repeller}  $A$-endomorphism $f\colon M^n\to M^n$ is an $A$-diffeomorphism, then the statement of the theorem is true only in case $n=1$. In case $n > 1$, basic set of a diffeomorphism of type $(n, 0)$ is a periodic source, thus it cannot be a submanifold of codimension one.
\end{rem}

\begin{rem}
Basic set $\Lambda$ from theorem~\ref{repeller} is not necessarily smoothly embedded in ambient manifold $M^n$. As an example we can consider an endomorphism of Riemann sphere induced by the map $z \to z^2 + c$. If parameter   $c$ is small enough but does not equal to zero, then  there is a basic set which is a repeller homeomorphic to the  circle but it is not smooth an any point (see, for example, \cite{Mi,Lu}).
\end{rem}

\begin{rem}
There exist repellers of type $(n, 0)$ that are not submanifolds.
It is well known that endomorphism of the Riemann sphere that are considered in holomorphic dynamics can have one-dimensional fractal repellers \cite{Mi,Lu}.
\end{rem}

\begin{teo}\label{attractor}
Let a  basic set $\Lambda$ be a codimension one topological submanifold of $M^n$. If $\Lambda$ is  an  attractor of an $A$-endomorphism $f\colon M^n\to M^n$ of  type $(n-1, 1)$, then:

1) $\Lambda$ is  smooth;

2) the restriction $f$ to $\Lambda$ is an expanding endomorphism.\footnote{In sense of definition~\ref{Shu_expand}.} 
\end{teo}

\begin{rem}
There exist attractors of type $(n-1, 1)$ that are not submanifolds. Let us consider an endomorphism of $3$-torus $f\colon \mathbb T^3 \to \mathbb T^3$ obtained as the direct product of an expanding endomorphism of the circle and a $DA$-diffeomorphism of the $2$-torus with one-dimensional attractor (see, for example, \cite{Sm}). Then nonwandering set of $f$ contains an attractor $\Lambda$ of type $(2, 1)$ locally homeomorphic to the product of the Cantor set and  the  $2$-dimensional disk.
\end{rem}

\begin{rem}
Examples of basic sets on $n$-manifolds defined in theorems~\ref{repeller} and~\ref{attractor} can be easily constructed. It is sufficient to consider a direct product of an expanding endomorphism of $(n-1)$-torus and a Morse-Smale diffeomorphism of the circle.
\end{rem}

\section{Auxiliary information}\label{aux}

Let $\langle \cdot, \cdot \rangle$ be a smooth Riemannian metric on manifold $TM^n$.  And $\rho$ be a metric on $M^n$ induced by $\langle \cdot, \cdot \rangle$.
Let $\Lambda$ be an invariant hyperbolic set of an endomorphism $f$. Then it is possible to introduce   local stable and unstable manifolds for points from $\Lambda$, as it was done for diffeomorphisms. However, there is a significant difference from the case of diffeomorphisms. For endomorphism local unstable manifold of a point $x \in \Lambda$ depend on $\bar x \in \hat x \cap \hat \Lambda$.

\begin{dfn}\label{def:local}
Let $\Lambda$ be hyperbolic invariant set of endomorphism $f \colon M^n \to M^n$ and $x \in \Lambda$, $\bar{x} \in \hat x \cap \hat \Lambda$. 
The set
$$W^s_{x, \varepsilon} = \{y \in M^n \mid \rho(f^n(x), f^n(y)) < \varepsilon, n = 0, 1, 2, \ldots\}$$
is called local stable manifold of point $x$, and the set
$$W^u_{x, \bar x, \varepsilon} = \{y \in M^n \mid \exists\, \bar y \in \hat y, 
\rho(x_n, y_n) < \varepsilon, n = 0, -1, -2, \ldots\}$$
is called local unstable manifold of a point $x$.
\end{dfn}

The structure of hyperbolic sets of endomorphisms was  studied in detail in  \cite{Prz1}. We will mention here some results important for the present paper.

\begin{prp}\label{metric}
For any hyperbolic set $\Lambda$ of endomorphism $f$ there exists a smooth Riemannian metric $\langle \cdot, \cdot \rangle_{\Lambda}$ on $TM^n$ equivalent to $\langle \cdot, \cdot \rangle$ and a real number $\lambda$, $0 < \lambda < 1$ such that for any $x \in \Lambda$ and $\bar x \in \hat \Lambda \cap \hat \Lambda$  following inequalities hold
$$\|Df_{x_i}(v)\|_{\Lambda} \leqslant \lambda \|v\|_{\Lambda}\quad\text{where} \quad v \in E^s_{x_i, \hat f^i(\bar x)},$$
$$\|Df_{x_i}(v)\|_{\Lambda} \geqslant (1/\lambda) \|v\|_{\Lambda}\quad\text{where} \quad v \in E^u_{x_i, \hat f^i(\bar x)},$$
$i \in \mathbb Z$.\end{prp}

\begin{prp}\label{local}
Let $\Lambda$ be a hyperbolic set of endomorphism $f$ then:
\begin{enumerate}
\item there exists $\varepsilon > 0$ such that for any $\bar x \in \hat x \cap \hat \Lambda$  the local stable $W^s_{x, \varepsilon}$ and the local unstable $W^u_{x, \bar x, \varepsilon}$ manifolds are smoothly embedded disks of topological dimension $\dim E^s_{x_0, \bar x}$ and $\dim E^u_{x_0, \bar x}$ tangent to $E^s_{x_0, \bar x}$ and $E^u_{x_0, \bar x}$ at point $x$;
\item $W^s_{x, \varepsilon}$ and $W^u_{x, \bar x, \varepsilon}$ depend continuously in $C^1$ topology on point $x$ and $\bar x$, respectively;
\item there exists $\mu > 1$ such that in metric $\rho$ on $M^n$ induced by Riemannian metric $\langle \cdot, \cdot \rangle_{\Lambda}$
\begin{enumerate}
\item for any points $y, z \in W^s_{x, \varepsilon}$ inequalities $\rho(f^{n+1}(y), f^{n+1}(z)) \leqslant (1/\mu) \rho(f^{n}(y), f^{n}(z))$, $n = 0, 1, 2, \ldots$ hold,
\item for any points $y, z \in W^u_{x, \bar x, \varepsilon}$ and $\bar y \in \hat y$, $\bar z \in \hat z$ satisfying inequalities from the definition of $W^u_{x, \bar x, \varepsilon}$, inequalities $\rho(y_{-n-1}, z_{-n-1}) \leqslant (1/\mu) \rho(y_{-n}, z_{-n})$, $n = 0, 1, 2, \ldots$ hold.
\end{enumerate}
\end{enumerate}
\end{prp}

\begin{prp}\label{convergence}
Let $f \colon M^n \to M^n$ be an $A$-endomorphism and $\Omega = \bigcup\limits_{j = 1}^l \Lambda_j$ its spectral decomposition. Then:
\begin{enumerate}
\item For any point $x\in M^n$ there exists unique basic set $\Lambda_{j_1}$ ($j_1 = \overline{1, l}$) such that $f^k(x) \to \Lambda_{j_1}$ as $k \to +\infty$. Moreover, there exists a point $y \in \Lambda_{j_1}$ such that $\rho(f^k(x), f^k(y)) \to 0$  as $k \to +\infty$.
\item For any $\bar x \in \hat M$ there exists a unique basic set $\Lambda_{j_2}$ ($j_1 = \overline{1, l}$) such that $x_i \to \Lambda_{j_2}$ as $i \to -\infty$. Moreover, there exists $\bar y \in \hat \Lambda_{j_2}$ such that $\rho(x_i, y_i) \to 0$ as $i \to -\infty$.
\end{enumerate}
\end{prp}

For $\delta > 0$ and for  $x \in M^n$  let $B_{\delta}(x)$ be an open ball of radius $\delta$ at point $x$ (i.e. $B_{\delta}(x) = \{y \in M^n \mid \rho(x, y) < \delta\}$) and $\bar B_{\delta}(x)$ be a closed ball of radius $\delta$ at point $x$ (i.e. $\bar B_{\delta}(x) = \{y \in M^n \mid \rho(x, y) \leqslant \delta\}$).

Next statement is a corollary of compactness. We present its prove here for the sake of completeness.

\begin{prp}\label{compact}
Let $K$ be a compact set of some metric space $X$ and $U$ be an open neighborhood of $K$. Then there exists $\delta > 0$ such that for any point $x\in K$ the inclusion $B_\delta(x) \subset U$ holds.
\end{prp}

\proof
Suppose that the statement of the lemma is not true. Then for any $\delta > 0$ there exists a point $x \in K$ such that $B_\delta(x) \setminus U \neq \varnothing$. Let us consider the sequence of positive numbers $\{ \delta_i \}_{i = 1}^{\infty}$  such that $\delta_i \to 0$ as $i \to \infty$ and a sequence of points $\{x_i\}_{i = 1}^{\infty}$ such that $B_{\delta_i}(x_i) \setminus U \neq \varnothing$.

Since $K$ is compact, without loss of generality we can assume that  the sequence $\{x_i\}_{i = 1}^{\infty}$ is convergent to some point $x_0 \in K$. Since $U$ is open, there exists  $\delta_0 > 0$ such that $B_{\delta_0}(x_0) \subset U$. Since sequences  $\{x_i\}_{i=1}^\infty$ and $\{\delta_i\}_{i=0}^\infty$ are convergent, there exists an integer $k \in \mathbb N$ such that inequalities $\rho(x_0, x_k) < \delta_0/2$ and $\delta_k < \delta_0/2$ hold.

Let us consider an arbitrary point  $y \in B_{\delta_k}(x_k)$. Then $\rho(x_0, y) \leqslant \rho(x_0, x_k) + \rho(x_k, y) < \delta_0/2 + \delta_k < \delta_0$  and $y \in B_{\delta_0}(x_0) \subset U$. It contradicts the assumption that $B_{\delta_k}(x_k) \setminus U \neq \varnothing$.
\qed

\section{The criterion for the existence of the attractor} 

Let $\Lambda$ be a basic set of $f$. Hereinafter we assume that $TM^n$ is equipped with a Riemannian metric $\langle \cdot, \cdot \rangle_\Lambda$ defined in proposition~\ref{metric} and that $\rho$ is a metric on $M^n$ induced by $\langle \cdot, \cdot \rangle_\Lambda$.

Theorem \ref{criterium} follows from lemmas \ref{attractor_manifolds},  \ref{k=n}  and  \ref{attractor_manifolds+}. 

\begin{lem}\label{attractor_manifolds}
Let $f$ be a $A$-endomorphism satisfying to  the axiom $ A $, $\Lambda $ be a basic set which  is an attractor. Then there exists $\varepsilon> 0$ such that for any point $x\in\Lambda$ and for any   $\bar x\in \hat x \cap \hat \Lambda$ the inclusion $W^u_ {x, \bar x, \varepsilon} \subset \Lambda$ holds.
\end{lem}
\proof 
Let $\Lambda$ be of type $(k, n-k)$. In case $k = 0$, lemma is trivial, since $W^u_ {x, \bar x, \varepsilon}$ coincides with a point $x$. That is why, we consider only the case $k \geqslant 1$.
Suppose that  $\Lambda $ does not coincide with the manifold  $M^n$ (otherwise the statement of the lemma  is obvious). Let $U$ be a neighborhood from  the definition of an attractor. Since $\Omega$ does not contain singular points (see item~1 of definition~\ref{axiom_a})  tangent map $Df_x \colon T_x M^n \to T_x M^n$ is nondegenerate at any point $x \in \Lambda$. So, it follows from inverse function theorem (see, for example, \cite{Zorich}, p. 499) that the restriction of $f$ to a sufficiently small neighborhood of the point $x$ is a local diffeomorphism. Therefore, for any $k\geqslant 0$, the set $f^k (U)$ contains an open neighborhood of the set  $\Lambda$.
Choose $\varepsilon >0$, which satisfies the conclusions of the proposition~\ref{local}. Suppose that the statement of the lemma is not true, then there is a point $x\in \Lambda $ and $\bar x \in \hat x \cap \hat \Lambda $ such that $W^u_{x, \bar x, \varepsilon} \not\subset \Lambda$. It follows from proposition~\ref{local} that $W^u_{x, \bar x, \varepsilon}$ is homeomorphic to $k$-disk. Therefore, there exists $y\in (U\setminus \Lambda) \cap W^u_{x, \bar x, \varepsilon}$. By the definition of attractor $f^l (y) \to \Lambda $ for $l \to + \infty$.
By definition of an attractor $\bigcap\limits_{j = 0}^{+\infty}f^{j}(U) = \Lambda$, therefore there exists $m\in \mathbb N$ such that $y\not\in f^m(U)$. It follows from item~3 of proposition~\ref{local} that for  the point $y$ there exists $\bar y \in \hat y$ such that $\rho (y_{l}, \Lambda) \to 0 $ as $l \to -\infty$.
This fact along with the fact that $f^m(U)$ contains an open neighborhood of $\Lambda $ imply an existence of a number $t \in \mathbb N $ such that $y_{-t} \in f^m(U) $ and $ f^{t} (y_{-t}) = y $.
Since $ f^t (f^m (U)) \subset f^m(U)$, it follows that $y = f^{t} (y_{-t}) \in f^{m} (U) $, which contradicts the choice of the number $ m $.
\qed

\begin{lem}\label{k=n}
Let $\Lambda$ be  basic set of an $A$-endomorphism of the type $(n, 0)$. If there exists $\varepsilon_1 > 0$ such that $W^u_{x, \bar x, \varepsilon_1} \subset \Lambda$ for any point $x \in \Lambda$ and for some $\bar x \in \hat x \cap \hat \Lambda$, then $\Lambda$ coincides with the ambient manifold $M^n$.
\end{lem}

\proof
Put $\varepsilon = \min\{\varepsilon_1, \varepsilon_2\}$, where $\varepsilon_2$ satisfies the conclusions of assertion~\ref{local}.
According to the  spectral decomposition theorem (see proposition  \ref{spectral}) $\Lambda$ is a closed set. Since $\Lambda$ is of type $(n, 0)$,  then  by  proposition~\ref{local} local unstable manifold $W^u_{x, \bar x, \varepsilon}$ contains an open $n$-dimensional disk. Therefore, $\Lambda$ is an open set. Thus, $\Lambda$ is simultaneously open and closed set,  hence it  coincides with $M^n$.
\qed

\begin{lem}\label{attractor_manifolds+}
Let $\Lambda$ be a basic set of an $A$-endomorphism of the type $(k, n-k)$, $0 \leqslant k \leqslant n - 1$. If there exists $\varepsilon_1 > 0$ such that $W^u_{x, \bar x, \varepsilon_1} \subset \Lambda$ for any point $x \in \Lambda$ and some $\bar x \in \hat x \cap \hat \Lambda$, then $\Lambda$ is an attractor.
\end{lem}

\proof
Put $\varepsilon = \min\{\varepsilon_1, \varepsilon_2\}$, where $\varepsilon_2$ satisfies conclusions of proposition~\ref{local}.
Let us show that there exists $\delta > 0$ such that for any point $x\in \Lambda$ and for any $\eta$ satisfying inequalities $0 <\eta \leqslant \delta$:
\begin{enumerate} 
\item\label{1} the intersection of $\bar B_{\eta}(x) \cap W^s_{\varepsilon}$ consists of one connected component,
\item\label{2} the intersection  $\partial \bar B_{\eta}(x) \cap W^s_{x, \varepsilon}$ is homeomorphic to $(n - k - 1)$-dimensional sphere.
\end{enumerate}

Assume the contrary. Since $W^s_{x, \varepsilon}$ is a smoothly embedded $(n - k)$-dimensional disk, for any point $x \in \Lambda$ there exists $\delta(x) > 0$ such that for any $\eta$ satisfying inequality $0 <\eta \leqslant \delta(x)$ properties \ref{1}, \ref{2} holds.

Consider a sequence $\{\delta_i\}_{i=1}^{+\infty}$ such that $\delta_i > 0$ and $\delta_i \to 0$ as  $i \to \infty$. The contrary assumption implies that for any $\delta_i$ there exists at least one point $x_i \in \Lambda$, such that $\delta_i$ does not satisfy at least one of the properties \ref{1}, \ref{2} at the point $x_i$. Since $\Lambda$ is compact, without loss of generality we can assume that  the sequence $\{x_i\}_{i = 1}^{\infty}$ is convergent to some point $x_0 \in \Lambda$. For a point $x_0$ there exists $\delta_0 > 0$ satisfying conditions \ref{1}, \ref{2}. By  continuous dependence $W^s_{x, \varepsilon}$ on a point $x \in \Lambda$ in $C^1$ topology there exists a neighborhood $V$ of the point $x_0$ such that $\delta_0 / 2$ satisfies to the properties \ref{1}, \ref{2} for any point $x \in V \cap \Lambda$. It contradicts the fact that $V \cap \Lambda$ contains points of the sequence $\{x_i\}_{i=1}^{+\infty}$ such that $\delta_i < \delta_0/2$.

Let $\delta$ satisfies the properties \ref{1}, \ref{2} for any point $x \in \Lambda$. Put $U = \bigcup\limits_{x \in \Lambda} \left( \bar B_\delta(x) \cap W^s_{x, \varepsilon} \right)$.

Let us  show that $U$ is a closed neighborhood of the attractor. According to the  choice of $\delta$ for any point $x \in \Lambda$ the set $B_\delta(x) \cap W^s_{x, \varepsilon}$ is homeomorphic to the closed $(n - k)$-dimensional disk. Consider an arbitrary point $x \in \Lambda$ and $\bar x \in \hat x \cap \hat \Lambda$ such that the inclusion $W^u_{x, \bar x, \varepsilon} \subset \Lambda$ holds. Since the unstable manifold $W^u_{x, \bar x, \varepsilon}$ is a smoothly embedded open $k$-dimensional disk, and the stable manifold $W^s_{x, \varepsilon}$ depends on $x$ continuously in $C^1$ topology the set $\bigcup\limits_{x \in W^u_{x, \bar x, \varepsilon}} \left( \bar B_\delta(x) \cap W^s_{x, \varepsilon} \right) \subset U$ is homeomorphic to the direct product of the open $k$-dimensional disk and  $(n - k)$-dimensional closed disk. Thus $U$ contains an open neighborhood of the set $\Lambda$.

Let us  show that $U$ is a closed set. Consider an arbitrary point $y \in \mathop{\rm{cl}}(U)$. There exists a sequence of points $\{y_i\}_{i = 1}^{+\infty}$, $y_i \in U$ converging to $y_0$. By construction, for any element of the sequence  $\{y_i\}_{i = 1}^{+\infty}$ one can find a point $x_i \in \Lambda$ such  that $y \in W^s_{x, \varepsilon}$ and $\rho(x_i, y_i) \leqslant \delta$. Since $\Lambda$ is compact, without loss of generality  the sequence $\{x_i\}_{i = 1}^{+\infty}$ can be considered as convergent to some point $x_0 \in \Lambda$. Since the metric $\rho$ is continuous (as a map from the direct product $M^n \times M^n$ to $\mathbb R$) and the sequence of pairs $(x_i, y_i)$ is convergent in $M^n \times M^n$, inequality $\rho(x_0, y) \leqslant \delta$ holds as $\rho(x_i, y_i) \leqslant \delta$ for any $i \in \mathbb N$.
By continuous dependence of $W^s_{x, \varepsilon}$ on  a point $x\in \Lambda$ in $C^1$ topology the point $y$ belongs to  $W^s_{x_0, \varepsilon}$. Thus, $y \in \bar B_\delta(x_0) \cap W^s_{x, \varepsilon}$ and  hence  $x \in U$.

Let us show that $f(U) \subset \mathop{\rm{int}}(U)$. Consider an arbitrary point $y \in U$. By construction, there exists a point  $x \in \Lambda$ such that  $x \in \bar B_\delta(x_0) \cap W^s_{x, \varepsilon}$. By item 3 of proposition~\ref{local} inequality $\rho(f(y), f(x)) \leqslant \delta / \mu$ holds for some $\mu > 1$. Thus, $f(y)$ belongs to the interior of $(n - k)$-dimensional closed disk $\bar B_\delta(f(x)) \cap W^s_{f(x), \varepsilon}$. By continuous dependence of local stable manifolds on   the point in $C^1$ topology the inclusion $y \in \mathop{\rm{int}}U$ holds.

Equality $\bigcap\limits_{n = 0}^\infty f^n(U) = \Lambda$ follows from the fact that for any $y \in U$ one has $y \in W^s_{x, \varepsilon}$ for some $x \in \Lambda$ and from item~3 of proposition~\ref{local}.

Thus, $U$ is the desired neighborhood from the definition of an attractor.
\qed
\section {The structure of attractors of type $ (N-1,1) $ and repellers of type $ (N, 0) $ that are topological manifolds of codimension one}

\subsection {The proof of the theorem \ref{repeller}}

\begin{lem} \label {non_zero_2}
If $\Lambda$ is a compact hyperbolic set of type $(n, 0)$  of an endomorphism $f$ and $\varepsilon> 0$ satisfies the conclusions of the proposition~\ref{local}, then there exists $\delta> 0$ such that for any $x \in \Lambda$,  $\bar x \in \hat x \cap \hat \Lambda$ the inclusion  $B_\delta (x) \subset (W^u_{x, \bar x, \varepsilon})$ holds.
\end{lem}
\proof
Since $f$ is continuous and $\Lambda$ is compact, the set $\hat \Lambda$ is a closed subset in $\widetilde M$.  Since  $M^n $ is compact, the set $\widetilde M$ is also compact. Thus, $\hat \Lambda$ is compact as a closed subset of a compact space.

Suppose that the statement of the lemma is not true, then there exists a sequence $\{\bar x^i\}_{i = 1}^{\infty}$ of  points in $\hat \Lambda$
and a sequence of numbers $\{\delta_i \} _{i = 1}^{\infty} $  such that $ B_{\delta_i} (x) \setminus W^u_{x_0^i, \bar x^ i, \varepsilon} \neq \varnothing$ and $ \delta_i \to 0 $ as $ i \to \infty $. Without loss of generality, we assume that the original sequence $ \{\bar x^i \}_{i = 0}^\infty$ is convergent to point $\bar x^0 \in \hat \Lambda$. By item~1 of proposition~\ref{local} local unstable manifold $W^u_{x_0^0, \bar x^0, \varepsilon}$ contains  an open disk $B_{\delta_0}(x_0^0)$. By item~2 of the same proposition there exists an integer $N$ such that for any $i > N$ one has $B_{\delta_0/2}(x) \subset W^u_{x_0^i, \bar x^ i, \varepsilon}$. It is contradicts the choice of the sequences $\{\bar x^i \}_{i = 0}^\infty$ and $\{\delta_i\}_{i = 0}^\infty$.
\qed

\begin{lem}\label{distance}
Let $K \subset M^n$ be compact and any point $x \in K$ is regular with respect to $f$.  Then there exist $\varepsilon > 0$ and an open neighborhood $V$ of $K$ such that if $x', x''\in V$ and $\rho(x', x'') < \varepsilon$, then $f(x') \neq f(x'')$.
\end{lem}
\proof
Since for any point $x \in K$ the tangent map $Df_x \colon T_xM^n \to T_{f(x)}M^n$ is nondegenerate, it follows from the inverse function theorem that there exists a neighborhood $U(x)$ of a point $x$ such that the restriction $f|_{U(x)} \colon U(x) \to f(U(x))$ is a diffeomorphism. Let us consider the open cover $U=\bigcup\limits_{x \in K} U(x)$ of the set $K$ with such neighborhoods. Since $M^n$ is  normal space  there exists an open neighborhood  $V$ of the set $K$ such that $U \supset \mathop{\text{cl}\,(V)} \supset V \supset K$ (see~for example~\cite{Eng}).

As the  manifold $M^n$ is compact, the set $\mathop{\text{cl}\,(V)}$ is also compact.  By Lebesgue's lemma\footnote{Let $X$ be a compact metric space, and $\{U\}$ be  its open cover. Then there exists a real number $\lambda > 0$ (called a Lebesgue's number) such that any subset of $X$ with a diameter less than $\lambda$ lies entirely in some element of the cover $\{U\}$(see~for example~\cite{Eng}). This Lemma also holds for any compact subset of metric space and its open covering.} there exists a real number $\lambda > 0$ such that for any pair of points $x', x'' \in \mathop{\text{cl}\,(V)}$ such that $\rho(x', x'') < \lambda$, there exists an element of the cover $U(x)$ that contains $x'$ and $x''$. If one set  $\varepsilon = \lambda$, then $V$ will be the desired neighborhood.
\qed

\begin{lem}\label{counterimages}
Let a basic set $\Lambda$ of an $A$-endomorphism $f$ be a topological submanifold of codimension one. Then there exists a neighborhood $U$ of  $\Lambda$ such that $f^{-1}(\Lambda) \cap U = \Lambda$.
\end{lem}

\proof
Suppose the contrary, then there exists a sequence $\{x_i\}_{i=0}^{+\infty}$, $x_i \in M^n \setminus \Lambda$, such that $\rho(x_i, \Lambda) \to 0$ as $i \to +\infty$ and $f(x_i) \in \Lambda$ for any $i \in \mathbb N$. By compactness of $M^n$ the sequence $\{x_i\}_{i=0}^{+\infty}$ can be considered as convergent to some point $x_0 \in \Lambda$.  By inverse function theorem there exists a neighborhood $V_1$ of the point $x_0$  such that the restriction of $f$ to  $V_1$ is a local diffeomorphisms. Since $\Lambda$ is a topological $(n-1)$-dimensional submanifold, there exists a neighborhood $V_2$ of the  point $x_0$ such that  $V_2 \cap \Lambda$ is homeomorphic to $(n-1)$-dimensional disk. Let $V$ be an open disk such that t$x \in V \subset V_1 \cap V_2$. The set $V \cap \Lambda$ is also homeomorphic to an open $(n-1)$-dimensional disk. Since $\Lambda$ is invariant, $f(V \cap \Lambda) \subset f(V) \cap \Lambda$.  Since $x_i \to x_0$ as $i \to +\infty$ , and $f(x_i) \in \Lambda$, by continuity $f$ inclusion  $f(x_i) \in f(V) \cap \Lambda$ holds for sufficiently big numbers   $i$. Thus, the restriction of $f$ to  $V$ is not injective map. It contradicts the fact that $f|_V$ is local diffeomorphisms.
\qed

\begin{lem}\label{expanding}
Let a basic set $\Lambda$ of an $A$-endomorphism $f$ be a topological submanifold of codimension one. If $\Lambda$ is of type $(n, 0)$,  then there exists a neighborhood $Q$ of $\Lambda$ such that $f$ is expanding on $Q$ in terms of definition~\ref{Kat_expand}.\footnote{$f \colon X \to X$ is called expanding on $A \subset X$, if there exist $\varepsilon > 0$ and $\mu > 1$ such that for any $x, y \in A$, $x \neq y$ one has $\rho(f(x), f(y)) > \mu \rho(x, y)$.}
\end{lem}

\proof
Let $\varepsilon_1$ satisfy conclusions of  proposition~\ref{local}. Let neighborhood $\Pi\supset  \Lambda$ and constant $\varepsilon_2 > 0$ satisfy the conclusion of lemma~\ref{distance}. Since $\Lambda$ is a basic set of an $A$-endomorphism, there exits $x \in \Lambda$ such that $\bigcup\limits_{i = 0}^\infty f^i(x)$ is dense in $\Lambda$. By proposition~\ref{compact} there  exists $\varepsilon_3>0$  and  $\bar x \in \hat x \cap \hat \Lambda$ such that  the inclusion $\bigcup\limits_{i \in \mathbb Z} W^u_{x_i, \hat f^i(\bar x), \varepsilon_3} \subset \Pi$ holds.

Set $\varepsilon = \min \{\varepsilon_1, \frac{\varepsilon_2}{2}, \varepsilon_3\}$ and $V = \bigcup\limits_{i \in \mathbb Z}\left[W^u_{x_i, \hat f^i(\bar x), \varepsilon} \cap f^{-1} \left( W^u_{x_{i+1}, \hat f^{i+1}(\bar x), \varepsilon_1} \right) \right]$.

Let us check  that the set $V$ is  a neighborhood of the basic set $\Lambda$. Since $\bigcup\limits_{i = 0}^{\infty}f^i(x)$ is dense in $\Lambda$, it is sufficient to show that there exist $\delta > 0$ such that for any $x_i$ one has $B_\delta(x_i) \subset W^u_{x_i, \hat f^i(\bar x), \varepsilon} \cap f^{-1} \left( W^u_{x_{i+1}, \hat f^{i+1}(\bar x), \varepsilon_1} \right)$. By lemma~\ref{non_zero_2} there exist $\delta_1 > 0$ such that $B_{\delta_1}(x_i) \subset W^u_{x_i, \hat f^i(\bar x), \varepsilon}$ for any $i \in \mathbb Z$. Since $\varepsilon \leqslant \varepsilon_1$ one has $B_{\delta_1}(x_{i+1}) \subset W^u_{x_{i+1}, \hat f^{i+1}(\bar x), \varepsilon_1}$ for any $i \in \mathbb Z$. 
Show that there exist $\delta_2 > 0$ such that $f\left( B_{\delta_2}(x_i) \right) \subset W^u_{x_{i+1}, \hat f^{i+1}(\bar x), \varepsilon_1}$. Suppose the contrary. Choose some sequence $\{\eta_j\}_{j=1}^{+\infty}$ such that $\eta_j \to 0$ as $j \to +\infty$. Then for any $j \in \mathbb N$ there exist points $x_{i_j}$ and $y_j$ such that $\rho(x_{i_j}, y_j) < \eta_j$ and $f(y_j) \not\in W^u_{x_{i+1}, \hat f^{i+1}(\bar x), \varepsilon_1}$. Since $B_{\delta_1}(x_{i+1}) \subset W^u_{x_{i+1}, \hat f^{i+1}(\bar x), \varepsilon}$, then $\rho(f(y_j), f(x_{i_j})) = \rho(f(y_j), x_{{i_j} + 1}) > \delta_1$. Since $M^n$ is compact without loss of generality one can consider the sequences  $x_{i_j}$ and $y_j$ to be convergent to some point $x_0$. Since $\rho(f(x_0), f(x_0))=0$, one has a contradiction with continuity of $f$. Then one can set $\delta = \min\{\delta_1, \delta_2\}$.

Since $M^n$ is normal and $\Lambda$ is closed, there exists an open set $Q$ such that $V \supset  \mathop{\text{cl}\,(Q)} \supset Q \supset \Lambda$.
Let us show that $Q$ is a desired neighborhood.

It follows from  Lebesgue's lemma and compactness of $\mathop{\text{cl}\,(Q)}$ that there exists $\lambda > 0$ such that for any pairs of points $y, z \in Q$ such that $\rho(y, z) < \lambda$  there exists $x_i$ such that $y, z \in W^u_{x_i, \hat f^i(\bar x), \varepsilon}$.

By definition of $Q$ one has $f\left(W^u_{x_i, \hat f^i(\bar x), \varepsilon}\right) \subset W^u_{x_{i+1}, \hat f^{i+1}(\bar x), \varepsilon_1}$ and $f(y), f(z) \in W^u_{x_{i+1}, \hat f^{i+1}(\bar x), \varepsilon_1}$. Set $y^+ = f(y), z^+ = f(z)$. By item 3 of proposition~\ref{local} there exist $\bar y^+ \in \hat y^+$ and $\bar z^+ \in \hat z^+$ such that $\rho(y^+_i, x_{i+1}) < \varepsilon$, $\rho(z^+_i, x_{i+1}) < \varepsilon$ for $i \leqslant 0$.

Show that $y = y^+_{-1}$ and $z = z^+_{-1}$. Let $y' = y^+_{-1}$ and $z' = y^+_{-1}$. Then one has $\rho(x_i, y') < \varepsilon$ and 	$\rho(x_i, z') < \varepsilon$. 
Let us  show that these inequalities cannot hold for any preimages  of the points $y^+$ and $z^+$ different from the points  $y$ and $z$. By  the choice of $\varepsilon$ we have $\rho(x_i, y) < \frac{\varepsilon_2} {2}$, $\rho(x_i, y') < \frac{\varepsilon_2}{2}$ and by the triangle inequality $\rho(y, y') \leqslant \rho(x_i, y) + \rho(x_i, y') < \varepsilon_2$. Then by Lemma~\ref{distance} the equality  $y = y'$ holds. Similarly the equality  $z = z'$ holds.

Hence  for any points $y, z \in Q$ such that $\rho(y, z) < \lambda$  the inequality $\rho(f(y), f(z)) > \mu \rho(y, z)$ holds for some $\mu > 1$. Thus  the restriction of $f$ to  $Q$ is expanding  endomorphism in the sense of definition~\ref{Kat_expand}.

\qed

\proof[Proof of the theorem~\ref{repeller}] 

The fact that $f|_{\Lambda}$ is an expanding endomorphism follows from the lemma~\ref{expanding}.
Show that the $\Lambda$ has a neighborhood $U$ defined in definition~\ref{def_repeller}.

Consider $x \in \Lambda$ such that $\bigcup\limits_{i=0}^\infty f^{i}(x)$ is dense in $\Lambda$. Choose $\varepsilon$ satisfying the conclusions of proposition~\ref{local}. Fix any $\bar x \in \hat x \cap \hat \Lambda$. By lemma~\ref{non_zero_2} the set $V_1 = \bigcup\limits_{i \in \mathbb Z}W^u_{x_i, \hat f^i(\bar x), \varepsilon}$ is an open cover of the set $\Lambda$. Let $V_2$ be a neighborhood of $\Lambda$ satisfying the conclusion of lemma~\ref{counterimages}. Since $M^n$ is normal, and basic sets $\Lambda_1, \ldots \Lambda_l$ are closed, then there exist disjoint open neighborhoods $Q_1, \ldots, Q_l$ of basic sets. Set $V_3 = Q_s$, where $Q_s$ is a neighborhood of $\Lambda$.

Since $\Lambda$ is compact, it follows from proposition~\ref{compact} that there exists $\delta_1 > 0$,  $\delta_2 > 0$ and $\delta_3 > 0$ such that $B_{\delta_1}(x_i) \subset V_1$, $B_{\delta_2}(x_i) \subset V_2$ and $B_{\delta_3}(x_i) \subset V_3$ for any $i \in \mathbb Z$. Set  $\delta = \min\{ \delta_1 / \mu, \delta_2, \delta_3 \}$, where $\mu > 1$ satisfies the conclusion of item~3 of proposition~\ref{local}. Set $U = \bigcup\limits_{i \in \mathbb Z}B_\delta(x_i)$. By construction, $U$  is also an open cover  of the set $\Lambda$ and $U \subset V_1 \cap V_2$. 
Let us  show that  $U$ is the desired neighborhood of the repeller $\Lambda$,  i.e. that the following conditions are satisfied: 

1) $f(U) \supset \mathop{\text{cl}\,(U)}$ and 

2) $\bigcap\limits_{k = 0}^{+\infty}f^{-k}(U) = \Lambda$.

1) Set $W = \bigcup\limits_{i \in \mathbb Z}B_{\mu \delta}(x_i)$, then $\mathop{\text{cl}\,(U)} \subset W \subset V_1$. Show that inclusion $W \subset f(U)$ holds. Indeed, by definition of $W$ for any point $y \in W$ there exists $i \in \mathbb Z$ such that $\rho(y, x_i) < \mu\delta$, and hence $y \in W^u_{x_i, \hat f^i(\bar x), \varepsilon}$. By item~3 of proposition~\ref{local} there exists a point $y' \in f^{-1}(y)$ such that $\rho(y', x_{i-1}) \leqslant (1/\mu) \rho(y, x_i) < (1/\mu)(\mu\delta) = \delta$. Thus, $y' \in B_\delta(x_{i-1}) \subset U$ and $y = f(y') \in f(U)$.

2) Assume the contrary, i.e. $\bigcap\limits_{k = 0}^{+\infty}f^{-k}(U) \neq \Lambda$. Then there exists a point $x \in \bigcap\limits_{k = 0}^{\infty}f^{-k}(U) \setminus \Lambda$. By definition of the point $x$ for any $k \geqslant 0$ one has $f^k(x) \in U$. By proposition~\ref{convergence} for any point $x \in M^n$ the sequence $\{f^i(x)\}_{i = 0}^{+\infty}$ converges to some basic set $\Lambda_j$. Since all points from $U \setminus \Lambda$ are wandering and $\mathop{\rm{cl}} U \cap \Lambda_j = \varnothing$ for any basic set $\Lambda_j \neq \Lambda$,  one has $\rho(f^i(x), \Lambda) \to 0$ as $i \to +\infty$ By proposition~\ref{convergence} there exists $y \in \Lambda$ such that $\rho(f^i(x), f^i(y)) \to 0$ as $i \to +\infty$ and, by lemma~\ref{counterimages} $f^i(x) \neq f^i(y)$. It contradicts the results of lemma~\ref{expanding}.

\qed

\subsection {The proof of the theorem \ref{attractor}}

\proof[Proof of the theorem~\ref{attractor}]

Let us show that $\Lambda$ is a smooth submanifold of the manifold $M^n$. Since 
$\Lambda$ is a topological submanifold,  for a  sufficiently small neighborhood $U(x)$ of any point $x \in \Lambda$ the set $\Lambda \cap U(x)$ is  homeomorphic to   $(n-1)$-dimensional disk. Choose $ \varepsilon> 0 $ such  that it satisfies the conclusions of the proposition~\ref{local}. Since the basic set $\Lambda$ is attractor  of type $(n-1, 1)$, by lemma~\ref{attractor_manifolds} any point $x \in \Lambda$ has a neighborhood (in  $\Lambda$) coinciding with a local unstable manifold $W^u_{x, \bar x, \varepsilon}$. As  $W^u_{x, \bar x, \varepsilon}$ is a smoothly embedded $ (n-1) $-dimensional disk, then  the basic set $\Lambda$ is a smooth submanifold.

Let us  show that the restriction of $ f_\Lambda $ is an expanding  map.
As the fiber $E^u_ {x_0, \bar x}$ is tangent to $ W^u_ {x, \bar x, \varepsilon} $ it does not depend on   $\bar x \in \hat x \cap \hat \Lambda$. Hence, $T \Lambda = \bigcup \limits_ {x \in \Lambda} E^u_ {x_0, \bar x}$. Then it follows from the hyperbolicity of the set $\Lambda $ and the definition \ref {Shu_expand} that the restriction $f |_ {\Lambda}$ is an expanding  map. \qed

\end{document}